\documentclass[12pt]{article}
\usepackage[cp1251]{inputenc}
\usepackage[T2A]{fontenc}
\usepackage[english]{babel}
\usepackage{amsthm}
\usepackage{amssymb}
\usepackage{amsmath, amsfonts}
\newtheorem{theorem}{Theorem}

\newtheorem{cor}{Corollary}
\begin{document}

\title{Sharp Estimates for Norms of Functions from Conjugate Classes in Metrics $C$ and $L$}

\author{V. F. Babenko\thanks{This paper was published in Russian in Studies in Modern Problems of Summation and Approximation
of Functions and Their Applications, Collect. Sci. Works, Dnepropetrovsk, 1973, pp. 3 -5}}
\date{}

\maketitle

\begin{abstract}
Sharp estimates for $C$ - and $L$ - norms of functions that are
conjugate with functions from the classes $W^rH_\omega$ of periodic
functions having prescribed concave  majorant of moduli of continuity, as
well as sharp estimates for the differences of such functions are
obtained.

\end{abstract}

Let $W^rH_\omega\; (r=0,1,...;\; W^0H_\omega =H_\omega)$ be the
class of $2\pi$-periodic, $r$ times differentiable functions $f(x)$ such that the
modulus of continuity of its $r$-th derivative $\omega(f^{(r)},t)\;
(f^{(0)}=f)$ is less or equal to the given modulus of continuity
$\omega (t)$. In addition, let $W_0^rH_\omega$ be the set of functions $f\in
W^rH_\omega$ satisfying
$$
\int\limits_0^{2\pi}f(t)dt=0.
$$
For a function $f(x)$, denote by $\widetilde{f}(x)$ the function that
is trigonometrically conjugate (see~\cite{Bary}) with $f(x)$. Let
$$
\widetilde{W}^rH_\omega :=\{ \widetilde{f}\; :\; f\in
{W}^rH_\omega\}.
$$

In papers~\cite{Korn1},~\cite{korn2} for any $r=0,1,...$ and concave
modulus of continuity $\omega (t)$ exact values
$$
M_r(\omega )_C:=\sup\limits_{f\in W_0^rH_\omega}\|
f\|_C,\;\;\Omega_r(\omega, t):=\sup\limits_{f\in
W^rH_\omega}\max\limits_x| f(x+t) - f(x)|
$$
and
$$
M_r(\omega
)_L:=\sup\limits_{f\in W_0^rH_\omega}\| f\|_L
$$
were obtained. Estimates for the values
$$
\widetilde{M}_r(\omega )_C:=\sup\limits_{f\in
\widetilde{W}_0^rH_\omega}\| f\|_C,\;\;\widetilde{\Omega}_r(\omega,
t):=\sup\limits_{f\in \widetilde{W}^rH_\omega}\max\limits_x| f(x+t)
- f(x)|, \; r\neq 0,
$$
as well as estimates for differences of higher order were found in~\cite{Motorn}.

We have proved the following statements.

\begin{theorem}
If $\omega(t)$ is concave modulus of continuity and for some
$\varepsilon > 0$ $$\int\limits_0^\varepsilon
\frac{\omega(t)}tdt<\infty,$$ then
\begin{enumerate}
\item
$$ \widetilde{M}_0(\omega)_C:=\sup\limits_{f\in H_\omega}\|
\widetilde{f}\|_C =\frac 1\pi\int\limits_0^{\pi /2}\frac{\omega
(2t)}{\sin t}dt,$$
\item
$$
\widetilde{\Omega}_0(\omega, t):=\sup\limits_{f\in
H_\omega}\max\limits_x| \widetilde{f}(x+t) - \widetilde{f}(x)|=\frac
2\pi\int\limits_0^{t/2}\frac{\sin t/2}{\cos u-\cos t/2}\omega(\rho
(u)-u)du,
$$
where $\rho (x)$ is defined by the equation
$$
\int\limits_0^{x}\frac{\sin t/2}{\cos u-\cos
t/2}du=-\int\limits_{\rho(x)}^{\pi}\frac{\sin t/2}{\cos u-\cos
t/2}du \;\;(0<x<t/2,\; t/2<\rho(x)<\pi).
$$
\end{enumerate}
\end{theorem}

\bigskip

Note that proof of this Theorem essentially uses the representation
(see~\cite[p. 573]{Bary})
$$
\widetilde{f}(x)=-\frac 1\pi\cdot
(P)\int\limits_{-\pi}^\pi\frac{f(t)}{2\tan \frac{t-x}{2}}dt
$$
and Lemma 5.1 from~\cite{Korn}.

\bigskip

\noindent {\bf Remark 1.} {\it There exist concave moduli of
continuity, e. g. $\omega(t)=-1/\ln t$, for which the
values $\widetilde{M}_0(\omega)_C$ and $\widetilde{\Omega}_0(\omega,
t)\;\; (t>0)$ are equal to $\infty$ .}

\bigskip

\noindent {\bf Remark 2.} {\it If $E_n(f)$ is the best uniform
approximation of $2\pi$ - periodic function $f$ by trigonometric
polynomials of order $n$, then
$$
\sup\limits_{f\in\widetilde{W}^rH_\omega}E_0(f)=\frac
12\max\limits_{t}\widetilde{\Omega}_r(\omega, t)=\frac
12\widetilde{\Omega}_r(\omega, \pi)\;\; (r=0,1,...).
$$}

\begin{theorem}
For any concave modulus of continuity the following equalities hold
$$
\widetilde{M}_r(\omega)_L:=\sup\limits_{f\in
\widetilde{W}^rH_\omega}\|
{f}\|_L=4\sum\limits_{i=0}^\infty\frac{(-1)^{i(r+1)}b_{2i+1}}{(2i+1)^{r+1}}
\;\; (r=1,2,...)
$$
where
$$
b_m=\frac 2\pi\int\limits_0^{\pi /2}\omega (2t)\sin mt\, dt.
$$
\end{theorem}

Proof of this Theorem uses analogs of Lemmas 1 and 2
from~\cite{korn2}.

\begin {cor}
For any $r=2,3,...$ the following equalities hold
$$
\sup\limits_{f\in {W}^rK}\| \widetilde{f}\|_L=\frac
{16K}\pi\sum\limits_{i=0}^\infty\frac{(-1)^{i(r+1)}}{(2i+1)^{r+2}}
$$
where $W^rK$ is the class of $2\pi$ - periodic functions $f(x)$ such
that $f^{(r-1)}$ is absolutely continuous and almost everywhere
$|f^{(r)}(x)|\le K$.

\end{cor}

\begin {cor}
For any $r=2,3,...$ and any concave modulus of continuity $\omega(t)$
the following equalities hold
$$
\sup\limits_{f\in
\widetilde{W}^rH_\omega}\bigvee\limits_0^{2\pi}(f)=
4\sum\limits_{i=0}^\infty\frac{(-1)^{ir}b_{2i+1}}{(2i+1)^{r}}.
$$
\end{cor}

\bigskip

\noindent {\bf Remark 3.} {\it For every presented above extremal
relation there exists a function from corresponding class that
realizes supremum}.

\end{document}